%% file: template-jems.tex
\documentclass{article}
\usepackage[journal=APDE,lang=british]{ems-journal-jems} 

\usepackage{textcomp}
\usepackage{mathrsfs}

\usepackage{chngcntr}
\usepackage{apptools}
\usepackage{verbatim}
\usepackage{tikz-cd}

\theoremstyle{definition}

\theoremstyle{remark}

\numberwithin{equation}{section}

\newcommand\restr[2]{{
		\left.\kern-\nulldelimiterspace 
		#1 
		\vphantom{\big|} 
		\right|_{#2} 
}}

\DeclareRobustCommand{\rchi}{{\mathpalette\irchi\relax}}
\newcommand{\irchi}[2]{\raisebox{\depth}{$#1\chi$}} 
\numberwithin{equation}{section}

\numberwithin{theorem}{section}
\numberwithin{proposition}{section}
\numberwithin{lemma}{section}
\numberwithin{remark}{section}
\numberwithin{exercise}{section}
\numberwithin{corollary}{section}

\begin{document}

\title{Nonlinear constrained optimization of Schur test functions}
\titlemark{Nonlinear constrained optimization of Schur test functions}


\emsauthor{1}{
	\givenname{Mikhail}
	\surname{Anikushin}
	\orcid{0000-0002-7781-8551}}{M.M.~Anikushin}
\emsauthor{2}{
	\givenname{Andrey}
	\surname{Romanov}
	\orcid{0009-0008-4886-2351}}{A.O.~Romanov}

\Emsaffil{1}{
	\pretext{}
	\department{Department of Applied Cybernetics}
	\organisation{Faculty of Mathematics and Mechanics, St Petersburg University}
	\rorid{023znxa73}
	\address{Universitetskiy prospekt 28}
	\zip{198504}
	\city{Peterhof}
	\country{Russia}
	\posttext{}
	\affemail{demolishka@gmail.com}
}
\Emsaffil{2}{
	\pretext{}
	\department{Department of Applied Cybernetics}
	\organisation{Faculty of Mathematics and Mechanics, St Petersburg University}
	\rorid{023znxa73}
	\address{Universitetskiy prospekt 28}
	\zip{198504}
	\city{Peterhof}
	\country{Russia}
	\posttext{}
	\affemail{romanov.andrey.twai@gmail.com}
}
	

	


\classification[90C26, 90C30, 90C90]{37L15}
\keywords{integral operators; Schur test; nonlinear constrained optimization; global stability}

\begin{abstract}
We apply the iterative nonlinear programming method, previously proposed in our earlier work, to optimize Schur test functions and thereby provide refined upper bounds for the norms of integral operators. As an illustration, we derive such bounds for transfer operators associated with twofold additive compound operators that arise in the study of delay equations. This is related to the verification of frequency inequalities that guarantee the global stability of nonlinear delay equations through the generalized Bendixson criterion.
\end{abstract}

\maketitle

\input{Introduction}
\input{Applications}
\input{Conclusion}





\begin{funding}
The reported study was funded by the Russian Science Foundation (Project 25-11-00147).
\end{funding}

\section*{Data availability}
The data that support the findings of this study can be generated using the scripts in the repository:\\
\centerline{\url{https://gitlab.com/romanov.andrey/freq-criterion-delay-schur-test}}

\section*{Conflict of interest}
The authors declare that they have no conflict of interest.


\end{document}

%% file: Introduction.tex
\section{Introduction}

Let $(\mathcal{X},\mu_{x})$ and $(\mathcal{Y},\mu_{y})$ be two measure spaces, and consider an integral operator $T_{K}$ with a nonnegative kernel $K=K(x,y)$. If there are positive functions $\mathfrak{p}$ and $\mathfrak{q}$ on $\mathcal{X}$ and $\mathcal{Y}$, respectively, and constants $\kappa_{x}, \kappa_{y} > 0$ such that
\begin{equation}
	\begin{split}
		&\int_{\mathcal{Y}}K(x,y) \mathfrak{q}(y)d\mu_{y}(y) \leq \kappa_{x} \mathfrak{p}(x) \qquad \text{for} \quad \text{$\mu_{x}$-almost all} \ x \in \mathcal{X},\\
		&\int_{\mathcal{X}}K(x,y) \mathfrak{p}(x)d\mu_{x}(x) \leq \kappa_{y} \mathfrak{q}(x) \qquad \text{for} \quad \text{$\mu_{y}$-almost all} \ y \in \mathcal{Y},
	\end{split}
\end{equation}
then the Schur test, see \cite[Theorem 5.2]{HalmosSunder1978}, asserts that $T_{K}$ is a well-defined bounded operator from $L_{2}(\mathcal{Y},\mu_{y})$ to $L_{2}(\mathcal{X},\mu_{x})$, and its norm is bounded from above by $\sqrt{\kappa_{x} \kappa_{y}}$.

Aimed to obtain refined upper bounds $\sqrt{\kappa_{x} \kappa_{y}}$ by considering the pair $(\mathfrak{p},\mathfrak{q})$ as a varying parameter, we obtain the minimax problem 
\begin{equation}
	\label{EQ: MinMaxSchurTest}
	\sup_{x \in \mathcal{X}, y \in \mathcal{Y}} \frac{\int_{\mathcal{Y}}K(x,y) \mathfrak{q}(y)d\mu_{y}(y)}{\mathfrak{p}(x)} \cdot \frac{\int_{\mathcal{X}}K(x,y) \mathfrak{p}(x)d\mu_{x}(x)}{\mathfrak{q}(y)} \to \min_{(\mathfrak{p},\mathfrak{q})}.
\end{equation}

In general, such kinds of optimization problems are not smooth and therefore demand the development of appropriate techniques, which are the topics of recent research \cite{CurtisOverton2012, CurtisMitchellOverton2017, KawanHafsteinGiesl2021}. However, in our work \cite{AnikushinRomanov2025RobustEstimates} we have proposed a simple general approach to such problems, which often reduces them, if not to smooth, then at least to more regular problems, which can be solved using standard methods of smooth nonlinear programming. We call this approach the \textit{iterative nonlinear programming}.

Namely, instead of resolving \eqref{EQ: MinMaxSchurTest}, we consider the analogous minimax problem, where the maximum is taken only over a finite set of \textit{reference points} $\{(x_{i}, y_{i})\}_{i=1}^{N}$ collected iteratively, see below. In this case, the problem can be reduced to an equivalent problem of nonlinear programming by means of introducing an auxiliary scalar parameter, whose value becomes the new objective function and bounds from above all the terms constituting the maximum, thereby defining $N$ inequality constraints. Then we complement such a problem by the requirement that the auxiliary parameter should change only slightly and resolve it via nonlinear programming techniques, among which we prefer the SLSQP method \cite{Kraft1988, JoshyHwang2024} implemented in the SciPy package for Python. This constitutes a single step of the iterative nonlinear programming. Here the simplest strategy for collecting reference points is to complement their set by the maximum found for the current parameter value at each iteration.

In practice, we are always restricted to multiparametric classes of functions. Suppose $\mathfrak{p}$ and $\mathfrak{q}$ are varied within such a class and depend smoothly on parameters. Then, for \eqref{EQ: MinMaxSchurTest}, it is immediately seen that we deal with smooth problems of nonlinear programming at each step of the optimization. For comparison, in \cite{AnikushinRomanov2025RobustEstimates}, the problem of optimizing singular values of derivatives of chaotic mappings over smooth families of metrics is considered. In this case, values at reference points may depend nonsmoothly on parameters. However, in concrete applications, their behavior is sufficiently regular, so the SLSQP method, although developed for smooth problems, succeeds in resolving the problem. Note also that the method\footnote{More precisely, its implementation in the PyGRANSO package for Python.} of \cite{CurtisMitchellOverton2017}, used instead of SLSQP, fails to converge\footnote{Since the method of \cite{CurtisMitchellOverton2017} does not have a well-established convergence theory, the problems considered in \cite{AnikushinRomanov2025RobustEstimates} may be a source of relevant counterexamples.} (and takes much time for this) for such kinds of problems\footnote{As for the method of \cite{CurtisOverton2012}, it is simply not suitable for solving such problems due to the need for many gradient evaluations (proportional to the number of parameters) caused by the gradient sampling.}. Similar observations about the qualitative and quantitative outperformance of SLSQP, when compared with known alternative techniques, are also reported in \cite{JoshyHwang2024}.

A rule of thumb is that for families with $r$ parameters, the iterative nonlinear programming optimization can be controlled using at most $O(r)$ reference points. Clearly, for effective control, their set must adapt to the current parameters. Depending on the problem, there may be different strategies for adapting reference points. For example, in the present paper, we simply discretize \eqref{EQ: MinMaxSchurTest} using uniform grids of points and apply the Simpson $1/3$-rule. Then, reference points are collected as maxima of the discretized expressions over the grid for the initial parameter values at each iteration. More complex strategies may involve nonuniform grids, whose nodes vary at each step, and adjustment of previously collected reference points to local maxima.

Given explicit formulas for the kernels and optimized Schur test functions, it should be possible to rigorously validate the upper estimates using interval arithmetic.

Our interest in the problem \eqref{EQ: MinMaxSchurTest} is related to transfer operators associated with twofold additive compound delay operators studied in our work \cite{AnikushinRomanov2023FreqConds}. Here, we deal with a family of operators depending on $\omega \in \mathbb{R}$, for which it is required to verify that their norms are strictly less than a certain threshold value $\Lambda^{-1}$, where $\Lambda > 0$. This represents one of the simplest frequency inequalities, which allows us to apply the generalized Bendixson criterion \cite{LiMuldowney1995LowBounds} that prevents closed invariant contours to exist on attractors of certain nonlinear scalar delay equations. Since such inequalities are preserved under small $C^{1}$-perturbations, the criterion can be also applied to $C^{1}$-close systems. From this, as in finite dimensions \cite{LiMuldowney1996SIAMGlobStab, Smith1986HD}, it is expected that these conditions imply the global stability, if appropriate variants of Pugh's closing lemma are developed\footnote{For our purposes, it is sufficient to obtain the closing via a $C^{1}$-small perturbation of the semiflow, which does not necessarily have to be generated by the same class of equations. For this reason, we strongly believe in the existence of such variants of the closing lemma in infinite dimensions.}.

Let us mention the theoretical basis behind our approach to the global stability of delay equations. It is constituted by the Frequency Theorem developed in \cite{Anikushin2020FreqDelay} (see also \cite{Anikushin2020FreqParab}) which guarantees the existence of certain quadratic Lyapunov-like functionals obtained through resolving appropriate infinite-horizon quadratic regulator problems, provided that the frequency inequalities are satisfied, and foundations for its applicability to compound cocycles generated by delay equations on exterior powers explored in \cite{Anikushin2023Comp}. In particular, such functionals can be used to investigate the uniform exponential stability of twofold compound cocycles. This is related to the problem of obtaining effective dimension estimates for delay equations, various aspects of which are discussed in our papers \cite{Anikushin2023LyapExp, AnikushinRomanov2023FreqConds,AnikushinRomanov2024EffEst}. For more general dynamical systems, see \cite{ZelikAttractors2022, KuzReit2020}.

In \cite{AnikushinRomanov2023FreqConds}, using finite-dimensional truncations of the integral operators, it is illustrated that the frequency-domain approach improves the global stability results for the Suarez--Schopf delayed oscillator, which were derived in \cite{Anikushin2022Semigroups} using the method of \cite{MalletParretNussbaum2013}, and for the Mackey--Glass equations, which can be obtained using the generalized Myshkis criterion from \cite{Lizetal2003}. In Section \ref{SEC: SchurTestDelayCompoundOperators}, we consider the latter example and justify the validity of frequency inequalities using upper estimates, which are more relevant to the problem.

Since the above approach is potentially applicable to a range of problems enjoying a kind of asymptotic compactness, such as certain parabolic, hyperbolic, or neutral delay equations, our results may be useful in some of these cases. However, we note that even in the case of scalar delay equations, the transfer operators associated with $m$-fold additive compound delay operators are not compact if $m > 2$. Therefore, the problem in general goes beyond the scope of integral operators with $L_{2}$-summable kernels. On the other hand, such integral operators still arise in the case of systems of delay equations and $m=2$, which is sufficient for studying the global stability. We plan to consider this case in future work.

More generally, the frequency-domain approach \cite{Anikushin2020FreqDelay, Anikushin2020FreqParab} is applicable for studying the existence of exponential dichotomies and inertial manifolds \cite{Anikushin2022Semigroups, AnikushinAADyn2021, AnikushinRom2023SS}, where similar frequency inequalities arise. In the case of inertial manifolds for delay equations, the inequalities are related to transfer operators acting between finite-dimensional spaces.

%% file: Applications.tex
\section{Norm estimates for transfer operators of additive compound delay operators}
\label{SEC: SchurTestDelayCompoundOperators}

In this section, we illustrate the method of iterative nonlinear programming to obtain refined upper bounds for the norm of integral operators associated with transfer operators of additive compound delay operators studied in our work \cite{AnikushinRomanov2023FreqConds}.

Let $a, b, \nu_{0} \in \mathbb{R}$, and $\tau > 0$ be given, and set $p = -\nu_{0} + i\omega$ for $\omega \in \mathbb{R}$. Consider the $2\times 2$ matrix $D_{0}$ along with its matrix exponential $e^{D_{0}t}$ for $t \in \mathbb{R}$ given by
\begin{equation}
	D_{0} = \begin{pmatrix}
		-a & b\\
		-be^{-p\tau} & a-p
	\end{pmatrix} \quad \text{and} \quad
	e^{D_{0}t} = \begin{pmatrix}
		g^{0}_{11}(t) & g^{0}_{12}(t)\\
		g^{0}_{21}(t) & g^{0}_{22}(t)
	\end{pmatrix}.
\end{equation}
It is well known that the entries of $e^{D_{0}t}$ can be computed explicitly according to the formula
\begin{equation}
	\label{EQ: ExplicitMatrixExponential}
	e^{D_{0} t} = e^{\alpha t} \left[ \left(\cosh(\delta t) - \alpha \frac{\sinh(\delta t)}{\delta}\right)I_{2} + \frac{\sinh(\delta t)}{\delta} D_{0} \right],
\end{equation}
where $\alpha \coloneq \operatorname{tr}D_{0}/2$, $\delta \coloneq \pm \sqrt{ - \det( D_{0} - \alpha I_{2}) }$, and $I_{2}$ denotes the identity $2\times 2$-matrix.

With the above quantities, we associate the kernel (depending on $\omega$)
\begin{equation}
	\label{EQ: ExampleResolventEquationsK2Formula}
	\begin{split}
		K(\theta,s) \coloneq \frac{e^{p\theta}e^{p\tau}g^{0}_{21}(\tau+\theta) (e^{-ps}g^{0}_{21}(-s) - g^{0}_{22}(\tau+s))}{1-e^{p\tau}g^{0}_{21}(\tau)} \\ 
		+\rchi_{[-\tau,\theta]}(s)e^{-ps}g^{0}_{21}(\theta-s) - \rchi_{[-\tau-\theta,0]}(s)g^{0}_{22}(\theta+\tau+s),
	\end{split}
\end{equation}
where $\theta, s \in [-\tau, 0]$, and $\rchi_{\mathcal{I}}$ denotes the characteristic function of the interval $\mathcal{I}$. It is correctly defined if $\nu_{0}$ is chosen such that the line $-\nu_{0} +i\mathbb{R}$ does not intersect the spectrum of the additive compound delay operator $A^{[\wedge 2]}$ studied in \cite[Section 4.5]{AnikushinRomanov2023FreqConds}.

From the asymptotic expansion of $\delta$ from \eqref{EQ: ExplicitMatrixExponential}, it can be shown that $K(\theta,s)$ is equivalent as $|\omega| \to \infty$ to the \textit{asymptotic kernel}
\begin{equation}
	\label{EQ: ExplicitComputationAsymptoticKernel}
	\bar{K}(\theta,s) \coloneq -e^{a \theta} \rchi_{[-\tau-\theta,0]}(s) e^{a (\tau+s)}e^{-p(\tau + s)}
\end{equation}
in the sense that $K(\theta,s) = \bar{K}(\theta,s) + O(|\omega|^{-1})$, where the remainder decays uniformly in $(\theta,s) \in [-\tau,0]^{2}$.

It can be seen that the norms of the integral operators $T_{\bar{K}}$ and $T_{|\bar{K}|}$ coincide. Clearly,
\begin{equation}
	\label{EQ: AsymptoticKernelModulusFormula}
	|\bar{K}|(\theta,s) = e^{a \theta} \rchi_{[-\tau-\theta,0]}(s) e^{(a+\nu_{0})(\tau+s)}
\end{equation}

Although there are explicit formulas for the kernels, it appears that even the norm of $T_{|\bar{K}|}$ cannot be explicitly computed\footnote{We refer to the \href{https://mathoverflow.net/q/500667}{discussion} on MathOverflow: https://mathoverflow.net/q/500667.}.

For most applications, computing the $L_{2}$-norm of $\bar{|K|}$ may be sufficient for validating the inequality $\|T_{K}\| < \Lambda^{-1}$ for all $\omega$ lying outside a sufficiently large segment $[-\Omega, \Omega]$, where $\Omega > 0$. Thus, the actual problem is to study what happens inside the segment.

Let $N_{1}$ and $N_{2}$ be real-valued functions on $[-\tau,0]$. Define the Schur test functions $\mathfrak{p}(\theta)$ and $\mathfrak{q}(s)$, where $\theta, s \in [-\tau,0]$, by
\begin{equation}
	\label{EQ: SchurTestFunctionsN1N2}
	\mathfrak{p}(\theta) \coloneq N_{1}(\theta)^{2} + 0.01 \quad \text{and} \quad \mathfrak{q}(s) \coloneq N_{2}(s)^{2} + 0.01.
\end{equation}

We consider a multiparametric family of $(N_{1},N_{2})$ determined by the neural network model with architecture $(1,30,2)$ and activation function $\sigma(y) = 1/(1+y^{2})$ for $y \in \mathbb{R}$. Specifically, as parameters we take $(30 \times 1)$- and $(2 \times 30)$-matrices $M_{1}$ and $M_{2}$, respectively, and $30$- and $2$-vectors $b_{1}$ and $b_{2}$, respectively. Then the model is given by
\begin{equation}
	\label{EQ: SchurTestNNmodel}
	N(x) = (N_{1}(x),N_{2}(x)) \coloneq M_{2}\sigma\left(M_{1}x + b_{1}\right) + b_{2},
\end{equation}
where the application of $\sigma$ to a vector is understood componentwise. Thus, there are $122$ parameters of the model.

We test the optimization algorithm by means of the parameters
\begin{equation}
	\label{EQ: ExplicitResolventParamatersForTest}
	a = -\tau' \gamma, \quad b=(\tau'\beta - \Lambda), \quad \tau=1, \quad \text{and} \quad \Lambda=\frac{1}{2}\tau'\beta\left( \frac{(\kappa-1)^{2}}{\kappa} +1\right),
\end{equation}
where $\gamma = 0.1$, $\beta = 0.2$, $\kappa = 10$, and $\tau' = 4.5$. These parameters correspond to the Mackey--Glass equations studied in \cite[Section 4.5]{AnikushinRomanov2023FreqConds}. We also take $\nu_{0} = 0.01$.

To discretize the problem for optimization, we use  a uniform grid of $251$ points partitioning the segment $[-\tau,0]$, and apply the Simpson $1/3$-rule.

\begin{figure}[t]
	\begin{minipage}{.5\textwidth}
		\includegraphics[width=\textwidth,angle=0]{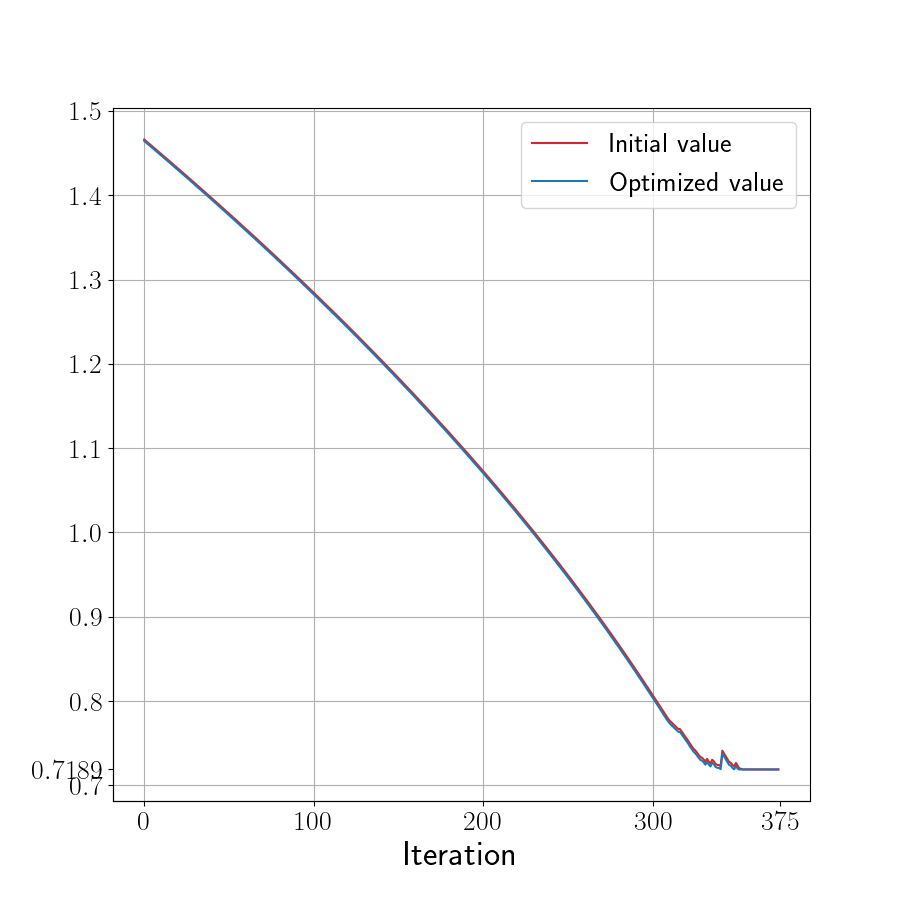}
	\end{minipage}%
	\begin{minipage}{.5\textwidth}
		\includegraphics[width=\textwidth,angle=0]{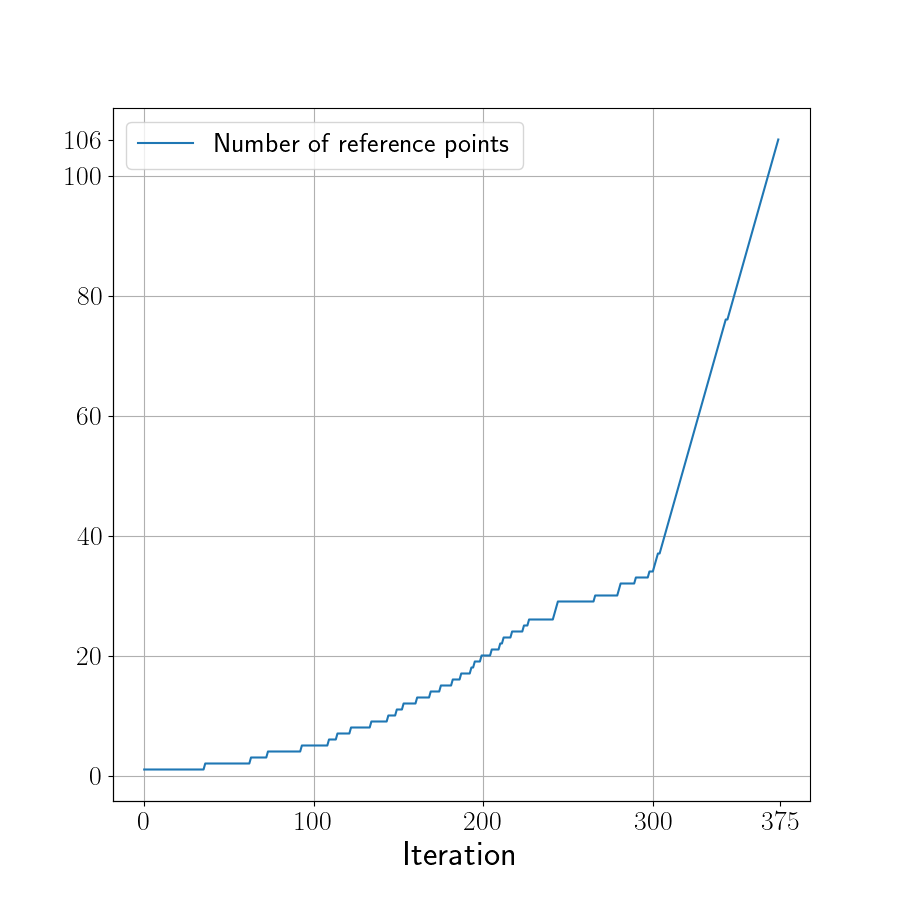}
	\end{minipage}%
	\caption{Dependence of the data versus iteration (horizontal) during the iterative nonlinear programming optimization of Schur test functions \eqref{EQ: SchurTestFunctionsN1N2}-\eqref{EQ: SchurTestNNmodel} applied to the kernel \eqref{EQ: ExampleResolventEquationsK2Formula} with parameters \eqref{EQ: ExplicitResolventParamatersForTest}, $\nu_{0} = 0.01$, and $\omega = 0$. (Left): Schur test estimate (red) and its optimized value over reference points (blue). (Right): Number of collected reference points (blue).}
	\label{FIG: CompoundMGTKernelShurTestOmegaZeroData}
\end{figure}

We first illustrate dynamics of the overall algorithm applied to estimate the norm\footnote{Note that we, in fact, estimate the norm of $T_{|K|}$.} of $T_{K}$ for $\omega = 0$. At each step of the iterative nonlinear programming, we bound the change of the objective function from below by $0.005$. Figure \ref{FIG: CompoundMGTKernelShurTestOmegaZeroData} shows the corresponding data obtained during the optimization, which successfully converged after 375 iterates, taking 35 seconds on our device\footnote{It is CPU AMD Ryzen 5 5600 OEM with 32GB 3200 MHz RAM.}. At the end, the algorithm collected 106 reference points, which is significantly smaller than the total number of points in the square grid, and comparable with the number of parameters.

Next, we apply the algorithm to obtain estimates for $\omega \in [-20, 20]$ with a step of $0.05$. When transiting from one value of $\omega$ to another, we keep the previous parameters of the neural network model, as well as the entire set of reference points, if there are less than $200$ of them, or leave only those which arose during the last $100$ iterations. This significantly lowers the number of iterations required to finish the optimization for the new value of $\omega$.

\begin{figure}[t]
	\begin{minipage}{.5\textwidth}
		\includegraphics[width=\textwidth,angle=0]{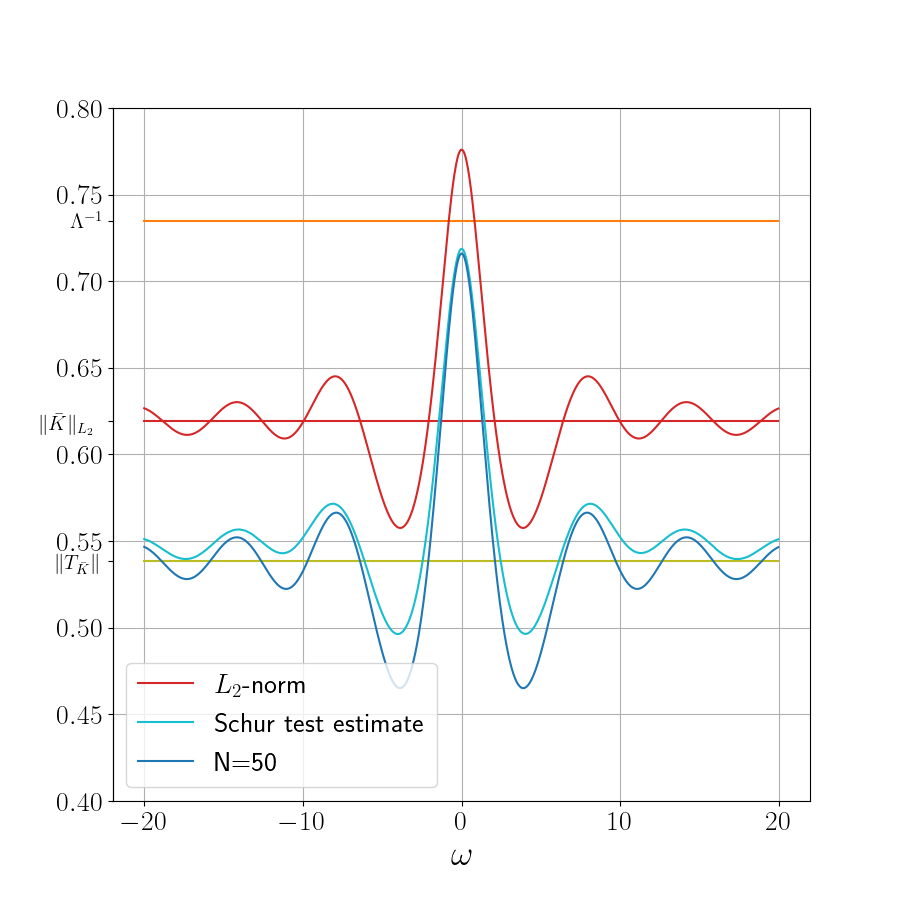}
	\end{minipage}%
	\begin{minipage}{.5\textwidth}
		\includegraphics[width=\textwidth,angle=0]{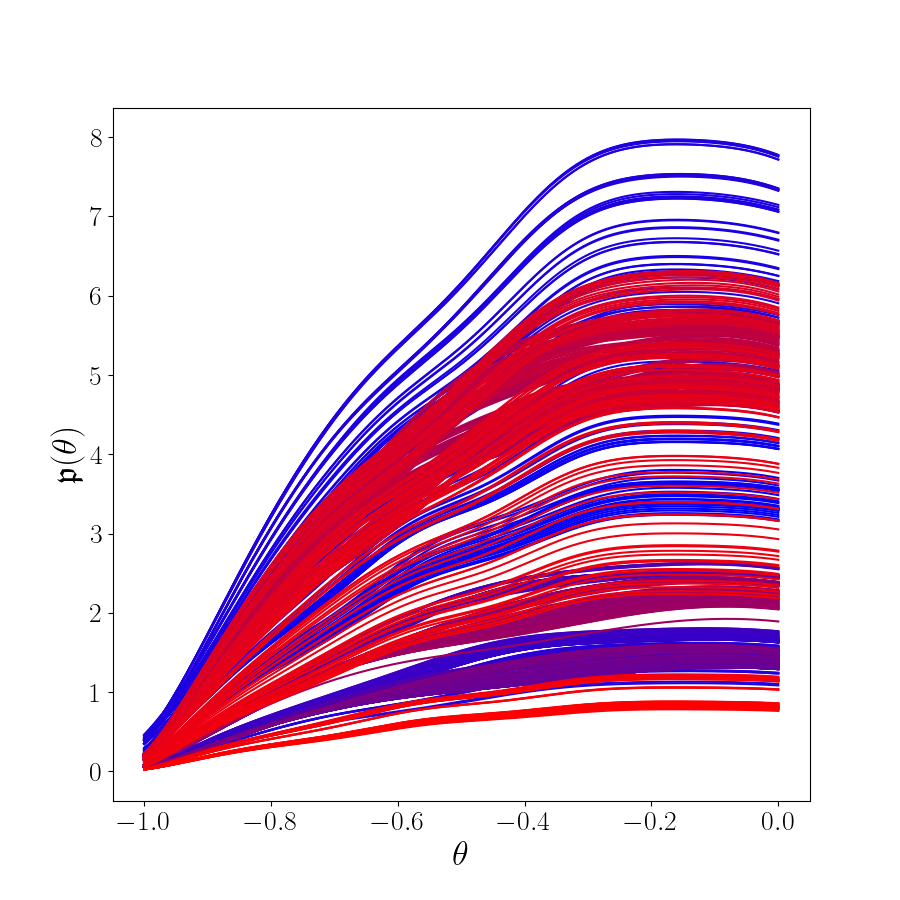}
	\end{minipage}%
	\caption{(Left): Estimates for the norm of $T_{K}$ versus $\omega$ (horizontal) in the case of \eqref{EQ: ExplicitResolventParamatersForTest} and $\nu_{0}=0.01$, namely, the $L_{2}$-norm of $K$ (red), the Schur test estimate (cyan), and the truncation norm (blue) for $N=50$ are shown. The horizontal lines pass through the threshold value $\Lambda^{-1}$ (orange), the $L_{2}$-norm of $\bar{K}$ (red), and the norm of $T_{\bar{K}}$ (olive) on the vertical axis. (Right): Plots of the optimized Schur test function $\mathfrak{p}$ for different $\omega$, colored using a blue-red colormap from $\omega = -20$ (blue) to $\omega = 20$ (red). See the repository for implementation details.}
	\label{FIG: CompoundMGTKernelApproximations}
\end{figure}

After optimization, we use the finer uniform grid of $1001$ points on $[-1,0]$ to test the result. We also compute the $L_{2}$-norm $\|K\|_{L_{2}}$ of the kernel $K$, which bounds from above the norm $\|T_{K}\|$ of $T_{K}$, and the norm of the truncation\footnote{This is the operator $P_{N} T_{K} P_{N}$, where $P_{N}$ is the orthogonal projector onto the subspace spanned by $\phi_{k}$ with $-N \leq k \leq N$.} of $T_{K}$ in the basis of trigonometric monomials $\phi_{k}(\theta)=e^{i2\pi k \theta}$ with $-N \leq k \leq N$ for $N=50$, which provides bounds from below for $\|T_{K}\|$. In their turn, both quantities $\|K\|_{L_{2}}$ and $\|T_{K}\|$ should be compared with their asymptotic values\footnote{For $\| \bar{K} \|_{L_{2}}$, there is an explicit formula that can be obtained from \eqref{EQ: AsymptoticKernelModulusFormula}. For $\|T_{\bar{K}}\|$, this is not the case, and here we use the truncation $P_{N}T_{K}P_{N}$ with $N=1000$. According to \eqref{EQ: ExplicitComputationAsymptoticKernel}, we can explicitly compute it, so no numerical integration is required.} $\| \bar{K} \|_{L_{2}}$ and $\|T_{\bar{K}}\|$, respectively. Figure \ref{FIG: CompoundMGTKernelApproximations} shows results of the optimization compared with these values. We point out that the Schur test estimate turns out to be very sharp in a neighborhood of $\omega=0$, where the most interesting behavior occurs.

%% file: Conclusion.tex
\section{Conclusion}
In this paper, we proposed the iterative nonlinear programming method for optimizing Schur test functions. By means of transfer operators associated with twofold compound delay operators, we demonstrated applications of the method to obtain refined norm estimates that can potentially be made rigorous using interval arithmetic. This provides a more reliable way for validating the associated frequency inequalities, which guarantee the global stability of certain scalar nonlinear delay equations through the generalized Bendixson criterion. This approach can also be applied to systems of delay equations, which we plan to consider in future work.